\documentclass[12pt]{report}
\usepackage{amssymb,amsmath,makeidx,verbatim}
\newcommand{\R}{\mathbb{R}}

\newcommand{\C}{\mathbb{C}}
\newcommand{\Z}{\mathbb{Z}}

\newcommand{\be}{\begin{enumerate}}
\newcommand{\ee}{\end{enumerate}}
\newcommand{\bq}{\begin{eqnarray*}}
\newcommand{\eq}{\end{eqnarray*}}

\begin{document}
\newcommand{\disp}{\displaystyle}
\thispagestyle{empty}
\begin{center}
\textsc{Non-spherical Harish-Chandra Fourier transforms on real reductive groups\\}
\ \\
\textsc{Olufemi O. Oyadare}\\
\ \\
Department of Mathematics,\\
Obafemi Awolowo University,\\
Ile-Ife, $220005,$ NIGERIA.\\
\text{E-mail: \textit{femi\_oya@yahoo.com}}\\
\end{center}
\begin{quote}
{\bf Abstract.} {\it The Harish-Chandra Fourier transform, $f\mapsto\mathcal{H}f,$ is a linear topological algebra isomorphism of the spherical (Schwartz) convolution algebra $\mathcal{C}^{p}(G//K)$ (where $K$ is a maximal compact subgroup of any arbitrarily chosen group $G$ in the Harish-Chandra class and $0<p\leq2$) onto the (Schwartz) multiplication algebra $\bar{\mathcal{Z}}({\mathfrak{F}}^{\epsilon})$ (of $\mathfrak{w}-$invariant members of $\mathcal{Z}({\mathfrak{F}}^{\epsilon}),$ with $\epsilon=(2/p)-1$). This is the well-known Trombi-Varadarajan theorem for spherical functions on the real reductive group, $G.$ Even though $\mathcal{C}^{p}(G//K)$ is a closed subalgebra of $\mathcal{C}^{p}(G),$ a similar theorem has not however been successfully proved for the full Schwartz convolution algebra $\mathcal{C}^{p}(G)$ except; for $\mathcal{C}^{p}(G/K)$ (whose method is essentially that of Trombi-Varadarajan, as shown by M. Eguchi); for few specific examples of groups (notably $G=SL(2,\R)$) and; for some notable values of $p$ (with restrictions on $G$ and/or on members of $\;\mathcal{C}^{p}(G)$). In this paper, we construct an appropriate image of the Harish-Chandra Fourier transform for the full Schwartz convolution algebra $\mathcal{C}^{p}(G),$ without any restriction on any of $G,p$ and members of $\;\mathcal{C}^{p}(G).$ Our proof, that the Harish-Chandra Fourier transform, $f\mapsto\mathcal{H}f,$ is a linear topological algebra isomorphism on $\mathcal{C}^{p}(G),$ equally shows that its image $\mathcal{C}^{p}(\widehat{G})$ can be nicely decomposed, that the full invariant harmonic analysis is available and implies that the definition of the Harish-Chandra Fourier transform may now be extended to include all $p-$tempered distributions on $G$ and to the zero-Schwartz spaces.}
\end{quote}
$\overline{2010\; \textmd{Mathematics}}$ Subject Classification: $43A85, \;\; 22E30, \;\; 22E46$

\ \\
{\bf \S1. Introduction.}

Let $G$ be a reductive group in the \textit{Harish-Chandra class} where $\mathcal{C}^{p}(G)$ is the \textit{Harish-Chandra-type} Schwartz algebra on $G,$ $0<p\leq2,$ with $\mathcal{C}^{2}(G)=:\mathcal{C}(G).$ It is known that $C^{\infty}_{c}(G)$ is dense in $\mathcal{C}^{p}(G),$ with continuous inclusion. The image of $\mathcal{C}^{p}(G)$ under the (Harish-Chandra) \textit{Fourier transform} on $G$ has been a pre-occupation of harmonic analysts since Harish-Chandra defined $\mathcal{C}(G)$ leading to the emergence of Arthur's thesis $[1a],$ where the Fourier image of $\mathcal{C}(G)$ was characterized for connected non-compact semisimple Lie groups of real rank one. Thereafter Eguchi $[3a.]$ removed the restriction of the real rank and considered non-compact real semisimple $G$ with only one conjugacy class of \textit{Cartan subgroups} as well as the Fourier image of $\mathcal{C}^{p}(G/K)$ in $[3b.],$ while Barker $[2.]$ considered $\mathcal{C}^{p}(SL(2,\R))$ as well as the zero-Schwartz space $\mathcal{C}^{0}(SL(2,\R)).$

The complete $p=2$ story for any real reductive $G$ is contained in Arthur $[1b,c].$ The most successful general result along the general case of $p$ is the well-known \textit{Trombi-Varadarajan Theorem} $[11.]$ which characterized the image of $\mathcal{C}^{p}(G//K),$ $0<p\leq2,$ for a maximal compact subgroup $K$ of a connected semisimple Lie group $G$ as a (Schwartz) multiplication algebra $\bar{\mathcal{Z}}({\mathfrak{F}}^{\epsilon})$ (of $\mathfrak{w}-$invariant members of $\mathcal{Z}({\mathfrak{F}}^{\epsilon}),$ with $\epsilon=(2/p)-1$); thus subsuming the works of Ehrepreis and Mautner $[5.]$ and Helgason $[7.].$ However the characterization of the image of $\mathcal{C}^{p}(G)$ for reductive groups $G$ in the Harish-Chandra class has not yet been achieved due to failure of the method of generalizing from the real rank one case (successfully employed in $[1b.,c.],$ $[3a.]$ and $[10c.]$) or from the spherical case (considered in $[11.]$).

This paper contains the full computation of the image of $\mathcal{C}^{p}(G)$ for reductive groups $G$ under the Harish-Chandra Fourier transform. It is organised as follows. The next section contains detailed preliminary matters concerning the structure of $G,$ its spherical functions and the Harish-Chandra-type Schwartz algebras, $\mathcal{C}^{p}(G).$ This section contains the most significant results on the spherical Harish-Chandra Fourier transform of these Schwartz algebras (Theorems $2.2$ and $2.3$) and also considered the system of differential equations satisfied by spherical functions, with a computation given for any real rank one $G.$

Our main results are contained in \S $3$ where a Schwartz algebra containing $\bar{\mathcal{Z}}({\mathfrak{F}}^{\epsilon})$ was constructed and we prove the \textit{full} non-spherical Harish-Chandra Fourier transforms of $\mathcal{C}^{p}(G)$ on real reductive groups $G$ (Theorem $3.8$) showing, at the same time, that its image $\mathcal{C}^{p}(\widehat{G})$ has nice \textit{decompositions} (Corollaries $3.9$ and $3.10;$ which confirm the real reason for the ease of transition of results from $[11.]$ to $[3b.]$) consisting of the \textit{Trombi-Varadrajan image,} $\bar{\mathcal{Z}}({\mathfrak{F}}^{\epsilon}),$ at its center. These decompositions save us the need to make endless asymptotic estimates in our analysis, showing that all such estimates have been subsumed in the Trombi-Varadarajan image $\bar{\mathcal{Z}}({\mathfrak{F}}^{\epsilon})$ (which enters the analysis naturally). The implication of this is that the fact that the spectrum in the theory for $G/K$ as computed in $[3b.]$ is still \textit{pure imaginary} is now shown to be mainly due to the contribution of the spherical case $G//K$ to the symmetric space $G/K$ (and not just \textit{carried over} to the case of $G/K,$ as posited in $[6.],\;p.\;355$). Indeed, the said decompositions of $\mathcal{C}^{p}(\widehat{G})$ give natural and direct paths to and from $\bar{\mathcal{Z}}({\mathfrak{F}}^{\epsilon}),$ as already evident from the results of $[3b.]$. We then show how the Trombi-Varadarajan theorem (Corollary $3.11$) could be recovered from our perspective. The Fourier transform of tempered distributions is thereafter extended to all of $\mathcal{C}^{p}(G)$ (Theorem $3.14$). We also lift the results of $[10d.]$ to give a full invariant harmonic analysis on $G$ (Theorems $3.18$ and $3.20$) with a proof of \textit{Rao-Varadarajan theorem} for $\mathcal{C}^{p}(G)$ (Theorem $3.21$). Basic results on the \textit{zero-Schwartz space} $\mathcal{C}^{0}(G)$ were considered (Theorem $3.22$) at the end of this section.

An application of our techniques is given in \S $4$ to (what we call) \textit{spherical convolutions,} $g_{\lambda,A},$ using the Harish-Chandra expansion of eigenfunctions on $G,$ thus leading to the Harish-Chandra Fourier transforms of a distinguished convolution subalgebra of $\mathcal{C}^{p}(G/K)$ (Theorem $4.3$) which contains the spherical part, $\mathcal{C}^{p}(G//K),$ of $\mathcal{C}^{p}(G).$ Further results on the structure of (canonical) wave-packets on $G$ and the considerations of the full \textit{Bochner theorem} shall be the subjects of future endeavours.

\ \\
{\bf \S2. Structure of the Schwartz algebras on $G.$}

Let $G$ be a group in the Harish-Chandra class. That is $G$ is a locally compact group with the properties that $G$ is reductive, with Lie algebra $\mathfrak{g},$ $[G:G^{0}]<\infty,$ where $G^{0}$ is the connected component of $G$ containing the identity, in which the analytic subgroup, $G_{1},$ of $G$ defined by $\mathfrak{g}_{1}=[\mathfrak{g},\mathfrak{g}]$ is closed in $G$ and of finite center and in which, if $G_{\C}$ is the adjoint group of $\mathfrak{g}_{\C},$ then $Ad(G)\subset G_{\C}.$ Such a group $G$ is endowed with a \textit{Cartan involution,} $\theta,$ whose fixed points form a \textit{maximal} compact subgroup, $K,$ of $G$ $[6.].$ $K$ meets all connected components of $G,$ in particular $K\cap G^{0}\neq\phi.$ Let $\mathfrak{t}$ denote the Lie algebra of $K.$

We denote the \textit{universal enveloping algebra} of $\mathfrak{g}_{\C}$ by $\mathcal{U}(\mathfrak{g}_{\C}),$ whose members may be viewed either as left or right invariant differential operators on $G.$ We shall write $f(x;a)$ for the left action $(af)(x)$ and $f(a;x)$ for the right action $(fa)(x)$ of $\mathcal{U}(\mathfrak{g}_{\C})$ on functions $f$ on $G.$ Let $\mathcal{C}(G)$ represents the space of $C^{\infty}-$functions $f$ on $G$ for which $$\sup_{x \in G}\mid f(b;x;a) \mid \Xi^{-1}(x)(1+\sigma(x))^{r}<\infty,$$ for $a,b \in \mathcal{U}(\mathfrak{g}_{\C})$ and $r>0.$ Here $\Xi$ and $\sigma$ are well-known \textit{elementary} spherical functions defined below on $G.$ $\mathcal{C}(G)$ is known to be a Schwartz algebra under convolution while $\mathcal{C}(G//K),$ consisting of the spherical members of $\mathcal{C}(G),$ is a closed commutative subalgebra. $C^{\infty}_{c}(G)$ is densely contained in $\mathcal{C}(G),$ with continuous inclusion.

Let $\hat{G}$ represent the set of equivalence classes of \textit{irreducible unitary representations} of $G.$ If $G_{1}$ is non-compact then the support of the Plancherel measure does not exhaust $\hat{G}.$ We write $\hat{G_{t}}$ for this support, which generally contains a discrete part, $\hat{G_{d}}$ ($\neq\emptyset,$ if $rank(G) = rank(K)$), and a continuous part, $\hat{G_{t}}\setminus \hat{G_{d}}$ ($\neq\emptyset,$ always).

If $\mathfrak{p}=\{X\in\mathfrak{g}:\theta X=-X\}$ then $\mathfrak{g}=\mathfrak{t}\oplus\mathfrak{p}.$ Choose a maximal abelian subspace  $\mathfrak{a}$ of $\mathfrak{p}$ with algebraic
dual $\mathfrak{a}^*$ and set $A =\exp \mathfrak{a}.$  For every $\lambda \in \mathfrak{a}^*$ put
$$\disp\mathfrak{g}_{\lambda} = \{X \in \mathfrak{g}: [H, X] =
\lambda(H)X, \forall  H \in \mathfrak{a}\},$$ and call $\lambda$ a \textit{restricted
root} of $(\mathfrak{g},\mathfrak{a})$ whenever $\mathfrak{g}_{\lambda}\neq\{0\}$.
Denote by $\mathfrak{a}'$ the open subset of $\mathfrak{a}$
where all restricted roots are $\neq 0$,  and call its connected
components the \textit{Weyl chambers}.  Let $\mathfrak{a}^+$ be one of the Weyl
chambers, define the restricted root $\lambda$ positive whenever it
is positive on $\mathfrak{a}^+$ and denote by $\triangle^+$ the set of all
restricted positive roots.  We then have the \textit{Iwasawa
decomposition} $G = KAN$, where $N$ is the analytic subgroup of $G$
corresponding to $\disp \mathfrak{n} = \sum_{\lambda \in \triangle^+} \mathfrak{g}_{\lambda},$
and the \textit{polar decomposition} $G = K\cdot
cl(A^+)\cdot K,$ with $A^+ = \exp \mathfrak{a}^+,$ and $cl(A^{+})$ denoting the closure of $A^{+}.$

If we set $\disp M = \{k
\in K: Ad(k)H = H, H\in \mathfrak{a}\}$ and $\disp M' = \{k
\in K : Ad(k)\mathfrak{a} \subset \mathfrak{a}\}$ and call them the
\textit{centralizer} and \textit{normalizer} of $\mathfrak{a}$ in $K,$ respectively, then;
(i) $M$ and $M'$ are compact and have the same Lie algebra and
(ii) the factor  $\mathfrak{w} = M'/M$ is a finite group called the \textit{Weyl
group}.  $\mathfrak{w}$ acts on $\mathfrak{a}^*_{\C}$ as a group of linear
transformations by the requirement $$(s\lambda)(H) =
\lambda(s^{-1}H),$$ $H \in \mathfrak{a}$, $s \in \mathfrak{w}$, $\lambda \in
\mathfrak{a}^*_\mathbb{\C}$, the complexification of $\mathfrak{a}^*$.  We then have the
\textit{Bruhat decomposition} $$\disp G = \bigsqcup_{s\in \mathfrak{w}} (B m_sB)$$ where
$B = MAN$ is a closed subgroup of $G$ and $m_s \in M'$ is the
representative of $s$ (i.e., $s = m_sM$).

Some of the most important functions on $G$ are the \textit{spherical
functions} which we now discuss as follows.  A non-zero continuous
function $\varphi$ on $G$ shall be called \textit{(elementary or zonal) spherical
function} whenever

$(i.)\;\varphi(e)=1,$ $$(ii.)\;\varphi \in C(G//K):=\{g\in
C(G): g(k_1 x k_2) = g(x), k_1,k_2 \in K, x \in G\}$$ and $(iii.)\;f*\varphi
= (f*\varphi)(e)\cdot \varphi$ for every $f \in C_c(G//K)$.  This
leads to the existence of a homomorphism $\lambda :
C_c(G//K)\rightarrow \C$ given as $\lambda(f) = (f*\varphi)(e)$.
This definition of an elementary spherical function is equivalent to the functional relation $$\disp
\int_K\varphi(xky)dk = \varphi(x)\varphi(y),$$ $x,y\in G$.  It has
been shown by Harish-Chandra [$6.$] that elementary spherical functions on $G$
can be parametrized by members of $\mathfrak{a}^*_{\C}$.  Indeed every elementary
spherical function on $G$ is of the form $$\disp
\varphi_{\lambda}(x) = \int_Ke^{(i\lambda-p)H(xk)}dk,\; \lambda
\in \mathfrak{a}^*_{\C},$$  $\disp \rho =
\frac{1}{2}\sum_{\lambda\in\triangle^+} m_{\lambda}\cdot\lambda,$ where
$m_{\lambda}=dim (\mathfrak{g}_\lambda),$ and that $\disp \varphi_{\lambda} =
\varphi_{\mu}$ iff $\lambda = s\mu$ for some $s \in \mathfrak{w}$.  Some of
the well-known properties are $\varphi_{-\lambda}(x^{-1}) =
\varphi_{\lambda}(x)$, $\varphi_{-\lambda}(x) =
\bar{\varphi}_{\bar{\lambda}}(x),\; \lambda \in \mathfrak{a}^*_{\C},\;\;
x \in G$, and if $\Omega$ is the \textit{Casimir operator} on $G$ then
$\Omega\varphi_{\lambda} = -(\langle\lambda,\lambda\rangle +
\langle \rho, \rho\rangle)\varphi_{\lambda},$ where $\lambda \in
\mathfrak{a}^*_{\C}$ and $\langle\lambda,\mu\rangle
:=tr(adH_{\lambda} \ adH_{\mu})$ for elements $H_{\lambda}$, $H_{\mu}
\in {\mathfrak{a}}.$ The elements $H_{\lambda}$, $H_{\mu}
\in {\mathfrak{a}}$  are uniquely defined by the requirement that $\lambda
(H)=tr(adH \ adH_{\lambda})$ and $\mu
(H)=tr(adH \ adH_{\mu})$ for every $H \in {\mathfrak{a}}$ ( [$6.$],
Propositions $3.1.4,\;3.2.1,\;3.2.2$ and Theorem $3.2.3$).  Clearly $\Omega\varphi_0 = 0.$

Let $$\varphi_0(x):= \int_{K}e^{(-\rho(H(xk)))}dk$$ be denoted
as $\Xi(x)$ and define $\sigma: G \rightarrow \C$ as
$\sigma(x) = \|X\|$ for every $x = k\exp X \in G,\;\; k \in K,\; X
\in \mathfrak{a},$ where $\|\cdot\|$ is a norm on the finite-dimensional
space $\mathfrak{a}.$ These two functions are zeroth elementary spherical functions on
$G$ and there exist numbers $c,d$ such that $1 \leq \Xi(a)
e^{\rho(\log a)} \leq c(1+\sigma(a))^d.$ Also there exists $r_0
> 0$ such that $c_0 =: \int_G\Xi(x)^2(1+\sigma(x))^{r_0}dx
< \infty$ ($[6.],$ p. $254$). For each
$0 \leq p \leq 2$ define ${\cal C}^p(G)$ to be the set consisting of
functions $f$ in $C^{\infty}(G)$ for which $$\mu^{p}_{g_1,
g_2;m}(f) :=\sup_G|f(g_1; x ; g_2)|\Xi (x)^{-2/p}(1+\sigma(x))^m <
\infty$$ where $g_1,g_2 \in \mathfrak{U}(\mathfrak{g}_{\C}),$ the \textit{universal
enveloping algebra} of $\mathfrak{g}_{\C},$ $m \in \Z^+, x \in G,$
$f(x;g_2) := \left.\frac{d}{dt}\right|_{t=0}f(x\cdot(\exp tg_2))$
and $f(g_1;x) :=\left.\frac{d}{dt}\right|_{t=0}f((\exp
tg_1)\cdot x).$

We call ${\cal C}^p(G)$ the Schwartz-type space on $G$
for each $0 < p \leq 2$ and note that ${\cal C}^2(G)$ is the earlier
Harish-Chandra space ${\cal C}(G)$ of rapidly decreasing functions on
$G.$ The inclusions $$C^{\infty}_{c}(G)\subset\bigcap_{0<p\leq 2}{\cal C}^p(G) \subset {\cal C}^p(G)
\subset L^p(G)$$ are continuous and with dense images. It also follows that
${\cal C}^p(G) \subseteq {\cal C}^q(G)$ whenever $0 \leq p \leq q
\leq 2.$ Each ${\cal C}^p(G)$ is closed under \textit{involution} and the
\textit{convolution}, $*.$ Indeed ${\cal C}^p(G)$ is a Fr$\acute{e}$chet algebra ($[12c.],\;p.\;357$) and the relation ${\cal C}^p(G)\ast{\cal C}^q(G)\subset{\cal C}^p(G)$ holds for all $p\geq q$ with $\frac{1}{p}+\frac{1}{q}=1;\;\;[3c.],$ Theorem $5.1.$ We endow ${\cal C}^p(G//K)$
with the relative topology as a subset of ${\cal C}^p(G).$

We shall say a function $f$ on $G$ satisfies a \textit{general strong inequality} if for any $r\geq0$ there is a constant $c=c_{r}>0$ such that
$$\mid f(y) \mid \leq c_{r} \Xi(y^{-1}x) (1+\sigma(y^{-1}x))^{-r}\;\;\;\;\;\forall\;x,y \in G.$$ We observe that if $x=e$ then, using the fact that $\Xi(y^{-1})=\Xi(y)$ and $\sigma(y^{-1})=\sigma(y),\;\forall\;y \in G,$ such a function satisfies $$\mid f(y) \mid \leq c_{r} \Xi(y^{-1}) (1+\sigma(y^{-1}))^{-r}=c_{r} \Xi(y) (1+\sigma(y))^{-r},\;\forall\;y \in G,$$ showing that a function on $G$ which satisfies a general strong inequality satisfies in particular a \textit{strong inequality} (in the classical sense of Harish-Chandra, $[12c.]$). Members of $\mathcal{C}(G)$ are those functions $f$ on $G$ for which $f(g_1; \cdot ; g_2)$ satisfies the strong inequality, for all $g_1,g_2 \in \mathfrak{U}(\mathfrak{g}_{\C}).$ We may then define $\mathcal{C}_{x}(G)$ to be those functions $f$ on $G$ for which $f(g_1; \cdot ; g_2)$ satisfies the general strong inequality, for all $g_1,g_2 \in \mathfrak{U}(\mathfrak{g}_{\C})$ and a fixed $x \in G.$ It is clear that $\mathcal{C}_{e}(G)=\mathcal{C}(G)$ and that $\bigcup_{x \in G}\mathcal{C}_{x}(G),$ which contains $\mathcal{C}(G),$ may be given an inductive limit topology.

\textbf{Proposition 2.1.} \textit{$\bigcup_{x \in G}\mathcal{C}_{x}(G)$ is a Schwartz algebra$.\Box$}

The algebra $\bigcup_{x \in G}\mathcal{C}_{x}(G)$ is worthy of an independent study. See $[9b.].$

For any measurable function $f$ on $G$ we define the \textit{Harish-Chandra Fourier
transform} $f\mapsto\mathcal{H}(f)$ as $\mathcal{H}(f)(\lambda) = \int_G f(x)
\varphi_{-\lambda}(x)dx,$ $\lambda \in \mathfrak{a}^*_{\C}=:\mathfrak{F}.$ We shall call it \textit{spherical} whenever $f\in\mathcal{C}^{p}(G//K)$ and, in this case, it may be shown that it is sufficient to define $\mathcal{H}f$ as $$(\mathcal{H}f)(\lambda):=\int_{AN}f(an)e^{(-\lambda+\rho)(\log a)}dadn,\;\lambda\in{\mathfrak{F}}_{I};$$ $[12a.],\;p.\;364.$

It is known (see $[6.]$) that for $f,g \in L^1(G)$ we have:
\begin{enumerate}
\item [$(i.)$] $\mathcal{H}(f*g) = \mathcal{H}(f)\cdot\mathcal{H}(g)$ on $ {\mathfrak{F}}^{1}$
whenever $f$ (or $g$) is right - (or left-) $K$-invariant; \item
[$(ii.)$] $\mathcal{H}(f^*)(\varphi) =
\overline{\mathcal{H}(f)(\varphi^*)},\;\varphi \in {\mathfrak{F}}^{1};$ hence
$\mathcal{H}(f^*) = \overline{\mathcal{H}(f)}$ on ${\cal P}:$ and, if we
define $f^{\#}(g) := \int_{K\times K}f(k_1xk_2)dk_1dk_2,  x\in
G,$ then \item [$(iii.)$] $\mathcal{H}(f^{\#})=\mathcal{H}(f)$ on ${\mathfrak{F}}^{1},$ where ${\mathfrak{F}}^{1}$ is the set of all bounded spherical functions and ${\cal P}$ is the subset of all positive-definite spherical functions.
\end{enumerate}

In order to know the image of the Harish-Chandra Fourier transform when
restricted to ${\cal C}^p(G//K)$ we need the following \textit{tube-spaces} that are central to the statement
of the well-known result of Trombi and Varadarajan [$11.$] (Theorem $2.2$ below).

Let $C_\rho$ be the closed convex hull of the (finite) set $\{s\rho :
s\in \mathfrak{w}\}$ in $\mathfrak{a}^*$, i.e., $$C_\rho =
\left\{\sum^n_{i=1}\lambda_i(s_i\rho) : \lambda_i \geq 0,\;\;\sum^n_{i=1}\lambda_i = 1,\;\;s_i \in \mathfrak{w}\right\}$$ where we recall that, for every
$H \in \mathfrak{a},$ $(s\rho)(H) = \frac{1}{2} \sum_{\lambda\in\triangle^+}
 m_{\lambda}\cdot\lambda (s^{-1}H).$  Now for each
$\epsilon > 0$ set ${\mathfrak{F}}^{\epsilon} = \mathfrak{a}^*+i\epsilon
C_\rho.$ Each ${\mathfrak{F}}^{\epsilon}$ is convex in $\mathfrak{a}^*_{\C}$ and
$$int({\mathfrak{F}}^{\epsilon}) =
\bigcup_{0<\epsilon'<\epsilon}{\mathfrak{F}}^{\epsilon^{'}}$$
([$11.$], Lemma $(3.2.2)$).  Let us define $\mathcal{Z}({\mathfrak{F}}^{0}) = \mathcal{S}
(\mathfrak{a}^*)$ and, for each $\epsilon>0,$ let
$\mathcal{Z}({\mathfrak{F}}^{\epsilon})$ be the space of all $\C$-valued
functions $\Phi$ such that  $(i.)$ $\Phi$ is defined and holomorphic
on $int({\mathfrak{F}}^{\epsilon}),$ and $(ii.)$ for each holomorphic
differential operator $D$ with polynomial coefficients we have $\sup_{int({\mathfrak{F}}^{\epsilon})}|D\Phi| < \infty.$ The space
$\mathcal{Z}({\mathfrak{F}}^{\epsilon})$ is converted to a Fr$\acute{e}$chet algebra by equipping it with the
topology generated by the collection, $\| \cdot \|_{\mathcal{Z}({\mathfrak{F}}^{\epsilon})},$ of seminorms given by $\|\Phi\|_{\mathcal{Z}({\mathfrak{F}}^{\epsilon})} := \sup_{int({\mathfrak{F}}^{\epsilon})}|D\Phi|.$  It is known that $D\Phi$ above extends to a continuous function on all of ${\mathfrak{F}}^{\epsilon}$
([$11.$], pp. $278-279$).  An appropriate subalgebra of
$\mathcal{Z}({\mathfrak{F}}^{\epsilon})$ for our purpose is the closed
subalgebra $\bar{\mathcal{Z}}({\mathfrak{F}}^{\epsilon})$ consisting of
$\mathfrak{w}$-invariant elements of $\mathcal{Z}({\mathfrak{F}}^{\epsilon}),$
$\epsilon \geq 0.$

{\bf 2.2 Theorem (Trombi-Varadarajan $[11.]$).}  \textit{Let $0 < p \leq 2$ and
set $\epsilon = \left(2/p\right)-1$.  Then the
Harish-Chandra Fourier transform $f \mapsto \mathcal{H}f$ is a linear
topological algebra isomorphism of ${\cal C}^p(G//K)$ onto $\bar{\mathcal{Z}}
({\mathfrak{F}}^{\epsilon}).\;\;\Box$}\\

For the Schwartz algebras ${\cal C}^p(G/K)$ a larger image than $\bar{\mathcal{Z}}({\mathfrak{F}}^{\epsilon})$ is required under the Harish-Chandra Fourier transform. Following Eguchi M. and Kowata A. $[4.]$ and Eguchi M. $[3b.]$ we define the space $\bar{\mathcal{Z}}(K/M\times{\mathfrak{F}}^{\epsilon})$ as the space of all $\mathfrak{w}$-invariant $C^{\infty}$ complex-valued functions $F$ on $K/M\times{\mathfrak{F}}_{I}$ which satisfy the following conditions:

$(i)$ for any $k\in K,$ the function $\lambda\mapsto F(kM:\lambda)$ extends holomorphically to $int({\mathfrak{F}}^{\epsilon});$

$(ii)$ for any $m\in\Z^{+},\;v \in S(\mathfrak{F}),$ $$\zeta^{\epsilon}_{v;m}(F) :=\sup_{(kM:\lambda)\in K/M\times int(\mathfrak{F})^{\epsilon}}|F(kM:\lambda ; \partial(v)|(1+ \mid \lambda\mid)^{m} <\infty.$$ The seminorms $\zeta^{\epsilon}_{v;m}$ restrict on $\bar{\mathcal{Z}}
({\mathfrak{F}}^{\epsilon})$ to the earlier Trombi-Varadarajan seminorms, $\| \cdot \|_{\mathcal{Z}({\mathfrak{F}}^{\epsilon})},$ and convert $\bar{\mathcal{Z}}(K/M\times{\mathfrak{F}}^{\epsilon})$ into a Fr$\acute{e}$chet space. Indeed, $\bar{\mathcal{Z}}({\mathfrak{F}}^{\epsilon})\subset\bar{\mathcal{Z}}(K/M\times{\mathfrak{F}}^{\epsilon}),$ as a closed subspace.

We define the map ${\cal C}^p(G/K)\rightarrow\bar{\mathcal{Z}}(K/M\times{\mathfrak{F}}^{\epsilon}):f\mapsto\mathcal{H}(f)$ now as $$\mathcal{H}(f)(kM:\lambda)=\int_{AN}f(kan)e^{(-\lambda+\rho)(\log a)}dadn,\;k\in K,\;\lambda\in\mathfrak{F}^{\epsilon},$$ referring to it as the \textit{symmetric} Harish-Chandra Fourier transform. A very important improvement on Theorem $2.2$ is the following.

{\bf 2.3 Theorem (Eguchi $[3b.]$).}  \textit{Let $0 < p \leq 2$ and
set $\epsilon = \left(2/p\right)-1$.  Then the
Harish-Chandra Fourier transform $f \mapsto \mathcal{H}f$ is a linear
topological algebra isomorphism of ${\cal C}^p(G/K)$ onto $\bar{\mathcal{Z}}(K/M\times{\mathfrak{F}}^{\epsilon}).\;\;\Box$}\\

Our first main result, given as Theorem $3.8,$ contains Theorems $2.2$ and $2.3$ as special cases.

The polar decomposition of $G$ implies that every $K-$biinvariant function on $G$ is completely determined by its restriction to $A^{+}.$ An example of such a function is the \textit{(zonal) spherical function,} $\varphi_{\lambda}, \lambda \in \mathfrak{a}^{\ast}_\mathbb{\C},$ on $G.$ If we denote the restriction of $\varphi_{\lambda}$ to $A^{+}$ as $\tilde \varphi_{\lambda},$ then the following system of differential equations hold: $$\tilde q \tilde \varphi_{\lambda}= \gamma (q)(\lambda) \tilde \varphi_{\lambda},$$ where $q \in \mathfrak{Q}(\mathfrak{g}_\mathbb{\C})(:= U(\mathfrak{g}_\mathbb{\C})^{K}$ = centralizer of K in $U(\mathfrak{g}_\mathbb{\C})$), $\gamma := \gamma_{\mathfrak{g}/\mathfrak{a}}$ is the \textit{Harish-Chandra homomorphism} of $\mathfrak{Q}(\mathfrak{g}_\mathbb{\C})$ onto $U(\mathfrak{g}_\mathbb{\C})^{\mathfrak{w}},$ the $\mathfrak{w}-$ invariant subspace of $U(\mathfrak{g}_\mathbb{\C}),$ with $\mathfrak{w}$ denoting the \textit{Weyl group} of the pair $(\mathfrak{g}, \mathfrak{a}),$ $\mathfrak{t} U(\mathfrak{g}_\mathbb{\C}) \bigcap \mathfrak{Q}(\mathfrak{g}_\mathbb{\C})$ is the kernel of $\gamma$ and $\tilde q$ is the restriction of $q$ to $A^{+}.$ Since $$ \widetilde{q \cdot f}
= \tilde q \cdot \tilde f,$$ for every $f \in C^{\infty}(G//K)$ we conclude that $\tilde q$ is the \textit{radial component} of $q.$ We define $q \in \mathfrak{Q}(\mathfrak{g}_\mathbb{\C})$ to be \textit{spherical} whenever $q = \tilde q.$

The above system of differential equations have been extensively used by Harish-Chandra in the investigation of the nature of the spherical functions, $\varphi_{\lambda},$ their asymptotic expansions and their contributions to the Schwartz algebras on $G.$ The history of this investigation dated back to the $1950's$ with the two-volume work of Harish-Chandra (See $[6.],\;p.\;190$), which still attracts the strength of twenty-first century mathematicians (See $[6.]$ and [$9a.$]). Other functions on $G$ satisfying different interesting transformations under members of $\mathfrak{Q}(\mathfrak{g}_\mathbb{\C})$ have also been studied in the light of the approach taken by Harish-Chandra. We refer to [$6.$] and the references cited in it for further discussion.

Now if $G$ is a semisimple Lie group with real rank $1$ then it is known (See $[11.]$) that the above system of differential equations can be replaced with $$\delta^{'}(\omega) \cdot \varphi_{\lambda}=\gamma (\omega)(\lambda) \cdot \varphi_{\lambda},$$ where $\omega$ is the \textit{Casimir operator} of $G$ and $\delta^{'}(\omega)$ denotes the radial component of the differential operator, $\delta^{'}(\omega),$ associated with $\omega.$ If we load the strucure of $G,$ as a real rank $1$ semisimple Lie group, into the last equation it becomes $$(\frac{d^{2}}{dt^{2}}+ \{(p+q) \coth t+ q \tanh t \} \frac{d}{dt}) f_{\lambda}=(\lambda^{2}- \frac{(p+2q)^{2}}{4})f_{\lambda},$$ where $p=n(\alpha),$ $q=n(2\alpha),$ $f_{\lambda}(t):= \varphi_{\lambda}(\exp tH_{0})$ and $H_{0}$ is chosen in $\mathfrak{a}$ such that $\alpha(H_{0})=1$ (See $[12b.],$ p. $190$ for the case of $G=SL(2, \R)$). Setting $z= -(\sinh t)^{2}$ transforms the above ordinary differential equation to the \textit{hypergeometric equation}$$(z(z-1) \frac{d^{2}}{dz^{2}} + ((a+b+1)z-c) \frac{d}{dt} + ab) g_{\lambda} = 0,$$ where $g_{\lambda}(z)=f_{\lambda}(t),$ $z<0,$ $a= \frac{p+2q+2\lambda}{4},$ $b= \frac{p+2q_2\lambda}{4}$ and $c= \frac{p+q+1}{2},$ whose solution is from here given by the well-known \textit{Gauss hypergeometric function,} $F(a,b,c:z),$ defined as $$F(a,b,c:z)= \sum^{\infty}_{k=0} \frac{(a)_{k} (b)_{k}}{(c)_{k}} \frac{z^{k}}{k!},$$ $\mid z \mid < 1.$ ([$6.$], p. $136$). It then follows that $$\varphi_{\lambda}(\exp tH_{0})= F(a,b,c:z)$$ (with $z=-(\sinh t)^{2}$), and we conclude that the spherical functions on real rank $1$ semisimple Lie groups are essentially the hypergeometric function. In other words, the hypergeometric functions form the spherical functions on any real rank $1$ semisimple Lie group.

In general and for any $G$ of arbitrary real rank we always have the Harish-Chandra series expansion for $\varphi_{\lambda}$ given as $$\varphi_{\lambda}(h)=\sum_{s\in\mathfrak{w}}c(s\lambda)\left(e^{(s\lambda-\rho)(\log h)}+\sum_{\mu\in L^{+}}a_{\mu}(s\lambda)e^{(s\lambda-\rho-\mu)(\log h)}\right),$$ valid for all $h\in A^{+},$ some $\lambda$ and $L^{+}$ as in $[6].$ The function $c$ is the well-known germ of harmonic analysis called Harish-Chandra $c-$function.
\ \\
\ \\
{\bf \S3. Harish-Chandra Fourier transform on $\mathcal{C}^{p}(G)$.}

We denote the set of equivalence classes of the necessarily finite-dimensional irreducible representations of $K$ by $\mathcal{E}(K)$ whose character is $\chi_{\mathfrak{d}},$ for every $\mathfrak{d}\in\mathcal{E}(K).$ The class functions $\xi_{\mathfrak{d}}:K\rightarrow\C$ defined as $\xi_{\mathfrak{d}}(k):=dim(\mathfrak{d})\chi_{\mathfrak{d}}(k^{-1})$ are idempotents (i.e., $\xi_{\mathfrak{d}}\ast\xi_{\mathfrak{d}}=\xi_{\mathfrak{d}}$ with $\xi_{\mathfrak{d}_{1}}\ast\xi_{\mathfrak{d}_{2}}=0$ whenever $\mathfrak{d}_{1}\neq\mathfrak{d}_{2}$). Choosing $\pi$ to be any representation of $K$ (which may be the restriction to $K$ of a representation of $G$) in a complete Hausdorff locally convex space, $V,$ a continuous projection operator on $V$ may be given as the image of $\xi_{\mathfrak{d}}$ under $\pi.$ That is, $$E_{\pi,\mathfrak{d}}:=\pi(\xi_{\mathfrak{d}})=
\int_{K}\xi_{\mathfrak{d}}(k)\pi(k)dk=dim(\mathfrak{d})\int_{K}\chi_{\mathfrak{d}}(k^{-1})\pi(k)dk$$ (Here $\int_{K}dk=1$) Idempotency of $\xi_{\mathfrak{d}}$ assures that $E_{\pi,\mathfrak{d}}$ is indeed a projection on $V$ (since $E^{2}_{\pi,\mathfrak{d}}=E_{\pi,\mathfrak{d}}$ and $E_{\pi,\mathfrak{d}_{1}}E_{\pi,\mathfrak{d}_{2}}=0$ whenever $\mathfrak{d}_{1}\neq\mathfrak{d}_{2}$) and that its range, written as $V_{\mathfrak{d}}(=E_{\pi,\mathfrak{d}}(V))$ is a closed linear subspace of $V$ consisting mainly of members of $V$ which \textit{transform according to} $\mathfrak{d};\;[12b.],\;p.\;109.$

The closed linear subspace $V_{\mathfrak{d}}$ above becomes familiar when $V=\mathcal{C}^{p}(G)$ under the usual regular representation. In this case the left and right regular representations are denoted as $l$ and $r$ given as $(l(x)f)(y)=f(x^{-1}y)$ and $(r(x)f)(y)=f(yx),$ respectively; $x,y\in G,\;f\in\mathcal{C}^{p}(G);$ and it may be computed that for any $\mathfrak{d}\in\mathcal{E}(K),$ $E_{l,\mathfrak{d}}=l(\xi_{\mathfrak{d}})$ and $E_{r,\mathfrak{d}}=r(\xi_{\mathfrak{d}})$ are the respective operators of left and right convolutions by the measure $\xi_{\mathfrak{d}}dk$ and $\overline{\xi}_{\mathfrak{d}}dk=\xi_{\overline{\mathfrak{d}}}dk,$ respectively. Here $\overline{\mathfrak{d}}$ is the class contragredient to $\mathfrak{d}.$ We therefore have a representation $l\times r$ of $G\times G$ on $\mathcal{C}(G)$ given as $((l\times r)(x,y)f)(z)=f(x^{-1}yz),$ $x,y,z\in G$ and the corresponding projection $E_{l\times r,(\mathfrak{d}_{1}\times\mathfrak{d}_{2})}=(l\times r)(\xi_{(\mathfrak{d}_{1}\times\mathfrak{d}_{2})}),$ which from the above remarks could be computed as $$E_{l\times r,(\mathfrak{d}_{1}\times\mathfrak{d}_{2})}f=\xi_{\mathfrak{d}_{1}}\ast f\ast\xi_{\mathfrak{d}_{2}},$$ $f\in\mathcal{C}^{p}(G).$

We now choose $\mathfrak{d}\in\mathcal{E}(K).$ The image of $\mathcal{C}^{p}(G)$ under $E_{l\times r,(\mathfrak{d}\times\overline{\mathfrak{d}})}$ is the closed subalgebra of $\mathcal{C}^{p}(G)$ denoted as $\mathcal{C}^{p}_{\mathfrak{d}}(G)$ and is exactly given as $$\mathcal{C}^{p}_{\mathfrak{d}}(G)=\xi_{\mathfrak{d}}\ast\mathcal{C}^{p}(G)\ast\xi_{\mathfrak{d}}=
\{\xi_{\mathfrak{d}}\ast f\ast\xi_{\mathfrak{d}}:\;f\in\mathcal{C}^{p}(G)\};$$ $[6.],\;p.\;11.$ Thus the members of $\mathcal{C}^{p}_{\mathfrak{d}}(G)$ are those of $\mathcal{C}^{p}(G)$ which may be written as $\xi_{\mathfrak{d}}\ast f\ast\xi_{\mathfrak{d}}$ for some $f\in\mathcal{C}^{p}(G).$ That is, every $g\in\mathcal{C}^{p}_{\mathfrak{d}}(G)$ is given as $g=\xi_{\mathfrak{d}}\ast f\ast\xi_{\mathfrak{d}},$  with $f\in\mathcal{C}^{p}(G).$ We shall henceforth write members $g$ of $\mathcal{C}^{p}_{\mathfrak{d}}(G)$ as $g_{\mathfrak{d},f},$ for some $f\in\mathcal{C}^{p}(G).$

\textbf{Lemma 3.1} \textit{Let $\mathfrak{d},\;\xi_{\mathfrak{d}}$ and $\mathcal{C}^{p}_{\mathfrak{d}}(G)$ be as above. Then every $g_{\mathfrak{d},f}\in\mathcal{C}^{p}_{\mathfrak{d}}(G)$ satisfies the transformation $\xi_{\mathfrak{d}}\ast g_{\mathfrak{d},f}\ast\xi_{\mathfrak{d}}=g_{\mathfrak{d},f}.$}

\textbf{Proof.} We know by definition that every $g_{\mathfrak{d},f}\in\mathcal{C}^{p}_{\mathfrak{d}}(G)$ is of the form $g_{\mathfrak{d},f}=\xi_{\mathfrak{d}}\ast f\ast\xi_{\mathfrak{d}};$  so that $\xi_{\mathfrak{d}}\ast g_{\mathfrak{d},f}\ast\xi_{\mathfrak{d}}=\xi_{\mathfrak{d}}\ast (\xi_{\mathfrak{d}}\ast f\ast\xi_{\mathfrak{d}})\ast\xi_{\mathfrak{d}}=\xi^{2}_{\mathfrak{d}}\ast f\ast\xi^{2}_{\mathfrak{d}}=\xi_{\mathfrak{d}}\ast f\ast\xi_{\mathfrak{d}}=g_{\mathfrak{d},f}.\;\Box$

It then means that members of the closed linear subspace $\mathcal{C}^{p}_{\mathfrak{d}}(G)$ are the $(\mathfrak{d},\overline{\mathfrak{d}})-$spherical functions in $\mathcal{C}^{p}(G),$ while the \textit{sphericalization operator} $E_{l\times r,(\mathfrak{d},\overline{\mathfrak{d}})}$ is the continuous projection $\mathcal{C}^{p}(G)\rightarrow\mathcal{C}^{p}_{\mathfrak{d}}(G).$ For the trivial representation $\mathfrak{d}=1$ of $K$ we shall write $\mathcal{C}^{p}_{\mathfrak{d}=1}(G)$ as $\mathcal{C}^{p}(G//K).$ Here the general transformation $\xi_{\mathfrak{d}}\ast g_{\mathfrak{d},f}\ast\xi_{\mathfrak{d}}=g_{\mathfrak{d},f}$ (for each $\mathfrak{d}\in\mathcal{E}(K),\;f\in\mathcal{C}^{p}(G)$), which is now $\xi_{1}\ast g_{1,f}\ast\xi_{1}=g_{1,f},$ becomes $(\xi_{1}\ast g_{1,f}\ast\xi_{1})(x)=g_{1,f}(x),$ $x\in G;$ leading to the the familiar expression $g_{1,f}(k_{1}xk_{2})=g_{1,f}(x),\;k_{1},\;k_{2}\in K,\;x\in G,$ for the $K-$biinvariance of spherical functions.

\textbf{Lemma 3.2} \textit{Let $0<p\leq 2$ and $f\in\mathcal{C}^{p}(G).$ Then $g_{1,f}\in\mathcal{C}^{p}(G//K).$}

\textbf{Proof.} The above remarks shows that $g_{1,f}$ is a spherical function on $G.$ If we now extend the definition of the function $\xi_{\mathfrak{d}}$ to all of $G$ by requiring that $\xi_{\mathfrak{d}}(kan)=e^{-(\lambda+\rho)(\log a)}\xi_{\mathfrak{d}}(k)$ with $\mathfrak{d}=1$ and note that $\mathcal{C}^{p}(G)$ is a convolution algebra (Theorem $5.1$ of $[3c.]$), we have the result$.\;\Box$

The last Lemma may be proved for the larger closed subalgebra $\mathcal{C}^{p}(G/K)$ of $\mathcal{C}^{p}(G)$ by the consideration of members of the closed subalgebra defined as $E_{r,\overline{\mathfrak{d}}}(\mathcal{C}^{p}(G))=\mathcal{C}^{p}(G)\ast\xi_{\mathfrak{d}}.$ The situation above may be completely extended to involve the idempotents $\xi_{F}$ defined  for any finite subset $F$ of $\mathcal{E}(K).$ In this case we set $\xi_{F}=\sum_{\mathfrak{d}\in F}\xi_{\mathfrak{d}}$ in order to have $E_{l\times r,F}$ and the closed linear subspace $\mathcal{C}^{p}_{F}(G)=\xi_{F}\ast\mathcal{C}^{p}(G)\ast\xi_{F}.$

The surjectivity of the map $E_{l\times r,(\mathfrak{d},\overline{\mathfrak{d}})}:\mathcal{C}^{p}(G)\rightarrow\mathcal{C}^{p}_{\mathfrak{d}}(G)$ assures that every Schwartz $(\mathfrak{d},\overline{\mathfrak{d}})-$spherical function on $G$ is in $\mathcal{C}^{p}_{\mathfrak{d}}(G).$ Hence, for any $f\in\mathcal{C}^{p}(G)$ the integral $\int_{G}g_{1,f}(x)\varphi_{\lambda}(x)dx$ converges absolutely and uniformly for all $\lambda\in\mathfrak{F}_{I},$ $([10c.],\;p.\;110,$ Lemma $3$)  and is continuous as a function on $\mathfrak{F}_{I},$ $([6.],\;p.\;262).$ Indeed by Theorem $2.2,$ the Harish-Chandra Fourier transform $\mathcal{H}:\mathcal{C}^{p}_{1}(G)\rightarrow\bar{\mathcal{Z}}
({\mathfrak{F}}^{\epsilon})$ is a linear topological algebra isomorphism ($[6.],\;p.\;354$) and for $f\in\mathcal{C}^{p}_{1}(G),$ the function $f\mapsto\mathcal{H}f$ is holomorphic on $int(\mathfrak{F}^{\epsilon}).$ This means that the inverses $(\mathcal{H}g_{1,f})^{-1}$ and $(\mathcal{H}\xi_{\mathfrak{d}})^{-1}$ exist as functions on (at least) $int(\mathfrak{F}^{\epsilon}).$

We actually have more than this, as contained in the following Lemma which gives an insight into the necessity of Eguchi's space $\bar{\mathcal{Z}}(K/M\times {\mathfrak{F}}^{\epsilon})$ over the Trombi-Varadarajan space $\bar{\mathcal{Z}}({\mathfrak{F}}^{\epsilon}),$ in the passage from $\mathcal{C}^{p}(G//K)$ (of Theorem $2.2$) to $\mathcal{C}^{p}(G/K)$ (of Theorem $2.3$). In what follows we shall denote the \textit{spherical} Harish-Chandra Fourier transform of $[11.],$ the \textit{symmetric} Harish-Chandra Fourier transform of $[3b.]$ and the \textit{general} Harish-Chandra Fourier transform of $[1c.]$ by the same symbol $\mathcal{H},$ since it will be clear which of the three is in use at any given place.

\textbf{Lemma 3.3.} \textit{The function $(\mathcal{H}\xi_{\mathfrak{d}})^{-1}$ exists and lives on $K/M\times int(\mathfrak{F}^{\epsilon}).$}

\textbf{Proof.} We first recall the extension of $\xi_{\mathfrak{d}}$ to all of $G (=KAN)$ by writing $\xi_{\mathfrak{d}}(kan)=e^{-(\lambda+\rho)(\log a)}\xi_{\mathfrak{d}}(k),$ so that $$(\mathcal{H}\xi_{\mathfrak{d}})(\lambda)=\int_{G}\xi_{\mathfrak{d}}(x)\varphi_{\lambda}(x)dx$$
$$=\int_{KAN}\xi_{\mathfrak{d}}(kan)\varphi_{\lambda}(kan)e^{2\rho(\log a)}dkdadn$$
$$=\int_{KAN}e^{-(\lambda+\rho)(\log a)}\xi_{\mathfrak{d}}(k)\varphi_{\lambda}(an)e^{2\rho(\log a)}dkdadn$$
$$=\int_{AN}e^{(-\lambda+\rho)(\log a)}\left[(\int_{K}\xi_{\mathfrak{d}}(k)dk)\varphi_{\lambda}(an)\right]dadn,$$ which by $(4.1.3)$ of $[3b.]$ is a function on $K/M\times int(\mathfrak{F}^{\epsilon}).$ We then have the result, due to its holomorphy on $int(\mathfrak{F}^{\epsilon}).\;\Box$

The situation of the last Lemma for $\mathfrak{d}=1$ is instructive and may be considered separately. Indeed, the same conclusion as in Lemma $3.3$ may be deduced for $\mathfrak{d}=1$ via the structure of the class-$1$ representations corresponding to the elementary spherical functions $\varphi_{\lambda}$ as follows.

\textbf{Lemma 3.4.} \textit{We always have that $(\mathcal{H}\xi_{1})(\lambda)=\widehat{(\mathcal{A}\varphi_{-\lambda})},$ where $\mathcal{A}$ is the Abel trasform and $\;\widehat{}\;$ is the Fourier transform on $A.$}

\textbf{Proof.} We already, know from Lemma $3.3,$ that $$(\mathcal{H}\xi_{1})(\lambda)=
\int_{AN}e^{(-\lambda+\rho)(\log a)}\varphi_{\lambda}(an)dadn.$$ Hence, $(\mathcal{H}\xi_{1})(\lambda)
=\int_{A}[e^{\rho(\log a)}\int_{N}\varphi_{\lambda}(an)dn]da
=\int_{A}(\mathcal{A}\varphi_{\lambda})(a)e^{-\lambda(\log a)}da=\widehat{(\mathcal{A}\varphi_{-\lambda})}.$ Thus the inverse $(\mathcal{H}\xi_{1})^{-1}(\lambda)$ is the inversion of the operation of finding the Abel transform of the elementary spherical functions $\varphi_{\lambda}$ (represented as an integral on $K/M$), followed by the operation of finding its Fourier transform on $A.$ Thus as $\lambda\in\mathfrak{F}_{I}$ is the parametrization of the class-$1$ representation $\pi_{1,\lambda}$

$(a.)$ whose matrix coefficients, defined by the constant function $1,$ is $\varphi_{\lambda}$ and
$(b.)$ which acts in $L^{2}(K/M)$ ($[6.],\;p.\;103 - 104$), it follows that $(\mathcal{H}\xi_{1})(\lambda)=(\pi_{1,\lambda}(\xi_{1})1,1)$ is in $L^{2}(K/M)$ and that the end result, $(\mathcal{H}\xi_{1})^{-1}(\lambda)$ of the inverted operations may be realized as a member of $L^{2}(K/M)$ for each $\lambda\in\mathfrak{F}_{I}.\;\Box$

The above two Lemmas reveal that the passage from the space $\bar{\mathcal{Z}}({\mathfrak{F}}^{\epsilon})$ of Trombi-Varadarajan to the space $\bar{\mathcal{Z}}(K/M\times {\mathfrak{F}}^{\epsilon})$ of Eguchi and the fact that the spectrum in $[3b.]$ is still \textit{pure imaginary} (as known in $[11.]$) are natural and are solely contributed to by the class-$1$ principal series representations, $\pi_{1,\lambda}.$ Due to the importance of these functions on $K/M$ (as seen in $[3b.]$) we shall now consider a candidate for the image of the whole of $\mathcal{C}^{p}(G)$ under the non-spherical Harish-Chandra Fourier transform, $\mathcal{H}.$

\textbf{Definition 3.5.} \textit{Let $0<p\leq 2$ and set $\mathcal{C}^{p}(\widehat{G})$ as $$\mathcal{C}^{p}(\widehat{G}):=
\{(\mathcal{H}\xi_{1})^{-1}\cdot h\cdot(\mathcal{H}\xi_{1})^{-1}:\;h\in\bar{\mathcal{Z}}({\mathfrak{F}}^{\epsilon})\}.\;\Box$$}
Observe that each of the three factors of every member of $\mathcal{C}^{p}(\widehat{G})$ is a function on $K/M\times {\mathfrak{F}}^{\epsilon},$ where we set $h(kM:\lambda):=h(\lambda),\;\lambda\in{\mathfrak{F}}^{\epsilon}.$ Thus every member of the above set $\mathcal{C}^{p}(\widehat{G})$ may therefore be seen as a function on $K/M\times {\mathfrak{F}}^{\epsilon}\times K/M$ ($[1c.],\;p.\;17$); so that $\mathcal{C}^{p}(\widehat{G})$ may be thought of as the function-space of the form $\bar{\mathcal{Z}}(K/M\times{\mathfrak{F}}^{\epsilon}\times K/M)$ (modelled on $\bar{\mathcal{Z}}({\mathfrak{F}}^{\epsilon})$ and $\bar{\mathcal{Z}}(K/M\times{\mathfrak{F}}^{\epsilon})$). It is clear, from Lemma $3.3,$ that $$\bar{\mathcal{Z}}({\mathfrak{F}}^{\epsilon})
\subset\bar{\mathcal{Z}}(K/M\times {\mathfrak{F}}^{\epsilon})\subset\mathcal{C}^{p}(\widehat{G})$$ and that both $\bar{\mathcal{Z}}({\mathfrak{F}}^{\epsilon})$ and $\bar{\mathcal{Z}}(K/M\times {\mathfrak{F}}^{\epsilon})$ are subspaces of $\mathcal{C}^{p}(\widehat{G}).$ We extend the seminorms on $\bar{\mathcal{Z}}(K/M\times {\mathfrak{F}}^{\epsilon})$ to all of $\mathcal{C}^{p}(\widehat{G})$ by setting $$\eta^{\epsilon}_{v;m}(F) :=\sup_{(kM:\lambda:kM)\in K/M\times int(\mathfrak{F})^{\epsilon}\times K/M}|F(kM:\lambda:kM ; \partial(v))|(1+ \mid \lambda\mid)^{m},$$ for $F\in\mathcal{C}^{p}(\widehat{G}),\;v \in S(\mathfrak{F}).$ Observe that $\eta^{\epsilon}_{v;m}$ restricts to $\zeta^{\epsilon}_{v;m}$ on $\bar{\mathcal{Z}}(K/M\times {\mathfrak{F}}^{\epsilon}).$ The following result is now immediate.

\textbf{Lemma 3.6.} \textit{$\mathcal{C}^{p}(\widehat{G})$ is a Schwartz algebra$.\;\Box$}

\textbf{Corollary 3.7.} \textit{$\bar{\mathcal{Z}}({\mathfrak{F}}^{\epsilon})$ and $\bar{\mathcal{Z}}(K/M\times {\mathfrak{F}}^{\epsilon})$ are closed subalgebras of $\mathcal{C}^{p}(\widehat{G}).$}

\textbf{Proof.} Both $\bar{\mathcal{Z}}({\mathfrak{F}}^{\epsilon})$ and $\bar{\mathcal{Z}}(K/M\times {\mathfrak{F}}^{\epsilon})$ are Schwartz subalgebras of $\mathcal{C}^{p}(\widehat{G}).\;\Box$

The following is our first main result giving the full non-spherical image of $\mathcal{C}^{p}(G)$ under the Harish-Chandra Fourier transform and may be compared with Theorems $2.2$ and $2.3.$

\textbf{Theorem 3.8.} (The Fundamental Theorem of Harmonic Analysis on $G$). \textit{Let $0<p\leq 2,$ then the Harish-Chandra Fourier transform, $\mathcal{H},$ sets up a linear topological algebra isomorphism $\mathcal{H}:\mathcal{C}^{p}(G)\rightarrow\mathcal{C}^{p}(\widehat{G}).$}

\textbf{Proof.} Let $f\in\mathcal{C}^{p}(G)$ be arbitrarily chosen, then (by Lemma $3.2$) there exists $g_{1,f}\in\mathcal{C}^{p}_{1}(G)$ given as $g_{1,f}=\xi_{1}\ast f\ast\xi_{1},$ such that (by Theorem $2.2$ we have) $\mathcal{H}g_{1,f}\in\bar{\mathcal{Z}}({\mathfrak{F}}^{\epsilon}).$ Hence, we have that $$\mathcal{H}(g_{1,f})=\mathcal{H}(\xi_{1}\ast f\ast\xi_{1})=\mathcal{H}(\xi_{1})\cdot\mathcal{H}(f)\cdot\mathcal{H}(\xi_{1}),$$ so that $$\mathcal{H}(f)
=(\mathcal{H}\xi_{1})^{-1}\cdot(\mathcal{H}g_{1,f})\cdot(\mathcal{H}\xi_{1})^{-1}\in\mathcal{C}^{p}(\widehat{G}).$$
That is, $\mathcal{H}(\mathcal{C}^{p}(G))\subset\mathcal{C}^{p}(\widehat{G}).$

Due to the detailed results of $[6.],\;[10c.],\;[3b.]$ and $[4.]$ on the already known properties of $\mathcal{H}$ concerning its linearity, continuity, injectivity and homomorphism, each of these properties of $\mathcal{H}$ now reduces to the same property for $\mathcal{H}_{|_{\mathcal{C}^{p}(G//K)}},$ due to the direct dependence of members of $\mathcal{C}^{p}(\widehat{G})$ on members of $\bar{\mathcal{Z}}({\mathfrak{F}}^{\epsilon})$ as now seen from Definition $3.5.$ For example, since $\mathcal{H}$ is already shown in the above paragraph to map $\mathcal{C}^{p}(G)$ to $\mathcal{C}^{p}(\widehat{G}),$ the injectivity of $\mathcal{H}$ follows from the already known injectivity of $\mathcal{H}_{|_{\mathcal{C}^{p}(G//K)}}.$ We are only left to show that $\mathcal{H}$ is surjective onto $\mathcal{C}^{p}(\widehat{G}).$

To this end, let $j\in\mathcal{C}^{p}(\widehat{G}).$ That is, let $j=(\mathcal{H}\xi_{1})^{-1}\cdot h\cdot(\mathcal{H}\xi_{1})^{-1}$ for some $h\in\bar{\mathcal{Z}}({\mathfrak{F}}^{\epsilon}).$ Define $f$ as the convolutions given as $$f=\mathcal{H}^{-1}((\mathcal{H}\xi_{1})^{-1})\ast\mathcal{H}^{-1}h\ast\mathcal{H}^{-1}((\mathcal{H}\xi_{1})^{-1}).$$
Now $\mathcal{H}^{-1}h\in\mathcal{C}^{p}(G//K)$ (by the Trombi-Varadarajan theorem); hence we have that $\mathcal{H}^{-1}h\in\mathcal{C}^{p}(G//K)\subset\mathcal{C}^{p}(G)$ and that $\mathcal{H}^{-1}((\mathcal{H}\xi_{1})^{-1})$ is the \textit{Eguchi-pullback} $\bar{\mathcal{Z}}(K/M\times {\mathfrak{F}}^{\epsilon})\rightarrow\mathcal{C}^{p}(G/K)\subset\mathcal{C}^{p}(G),\;[3b.],\;p.\;193.$

Hence, as $f$ is now shown to be the convolutions of members of $\mathcal{C}^{p}(G)$ and $\mathcal{C}^{p}(G)$ is a convolution algebra, we conclude that $$f:=\mathcal{H}^{-1}((\mathcal{H}\xi_{1})^{-1})\ast\mathcal{H}^{-1}h
\ast\mathcal{H}^{-1}((\mathcal{H}\xi_{1})^{-1})\in\mathcal{C}^{p}(G).$$ Finally, we have that $\mathcal{H}(f)=(\mathcal{H}\mathcal{H}^{-1}((\mathcal{H}\xi_{1})^{-1}))\cdot(\mathcal{H}\mathcal{H}^{-1}h)
\cdot(\mathcal{H}\mathcal{H}^{-1}((\mathcal{H}\xi_{1})^{-1}))$ $=(\mathcal{H}\xi_{1})^{-1}\cdot h\cdot(\mathcal{H}\xi_{1})^{-1}=j$ as expected. The inverse map $$\mathcal{H}^{-1}:\mathcal{C}^{p}(\widehat{G})\rightarrow\mathcal{C}^{p}(G)$$ is continuous, being the continuous extension to $\mathcal{C}^{p}(\widehat{G})$ of the (continuous) map $\mathcal{H}^{-1}_{|_{\bar{\mathcal{Z}}(K/M\times {\mathfrak{F}}^{\epsilon})}}:\bar{\mathcal{Z}}(K/M\times {\mathfrak{F}}^{\epsilon})\rightarrow\mathcal{C}^{p}(G/K)$ of Eguchi, $[3b.].\;\Box$

\textbf{Corollary 3.9.} \textit{The algebra $\mathcal{H}(\mathcal{C}^{p}(G))\cong\mathcal{C}^{p}(\widehat{G})$ may also be seen as $$(\mathcal{H}\xi_{1})^{-1}\cdot\bar{\mathcal{Z}}
({\mathfrak{F}}^{\epsilon})\cdot(\mathcal{H}\xi_{1})^{-1},$$ where $(\mathcal{H}\xi_{1})^{-1}$ denotes the fixed function given as $(\mathcal{H}\xi_{1})^{-1}(\lambda),$ for $\lambda\in{\mathfrak{F}}^{\epsilon}.$}

\textbf{Proof.} By Definition $3.5$ and Theorem $3.8$ we have that $\mathcal{H}(\mathcal{C}^{p}(G))\cong\mathcal{C}^{p}(\widehat{G})
:=\{(\mathcal{H}\xi_{1})^{-1}\cdot h\cdot(\mathcal{H}\xi_{1})^{-1}:\;h\in\bar{\mathcal{Z}}({\mathfrak{F}}^{\epsilon})\}
=(\mathcal{H}\xi_{1})^{-1}\cdot\bar{\mathcal{Z}}
({\mathfrak{F}}^{\epsilon})\cdot(\mathcal{H}\xi_{1})^{-1}.\;\Box$

Direct computations of members of $\mathcal{C}^{p}(\widehat{G})$ may also be embarked on. It is clear that the algebra $\mathcal{C}^{p}(\widehat{G})$ is still in its \textit{bundled} form and that it would need to be further opened up than has been done in the decomposition contained in Corollary $3.9.$ Indeed, it has to be explicitly computed and understood for different examples of $G$ and its members parametrized, while already known cases are shown to be deduced from it. In particular, we are yet to give an explicit computation of the inverse map $\mathcal{H}^{-1}:\mathcal{C}^{p}(\widehat{G})\rightarrow\mathcal{C}^{p}(G)$ in the theory of \textit{wave-packets} or discuss the explicit nature of contributions of the \textit{discrete} and \textit{principal series} of representations of $G$ to $\mathcal{C}^{p}(\widehat{G})$ or ask if there is still a split (as known for $\mathcal{C}^{p}(\widehat{SL(2,\R)})$ in $[2.]$ and for $\mathcal{C}^{p}(\widehat{G}:F)$ in $[10d.]$) into the discrete and principal parts at the level of $\mathcal{C}^{p}(\widehat{G})$ or consider other well known results on the \textit{spherical} case for all of $\mathcal{C}^{p}(G).$

However, Theorem $3.8$ marks a significant attainment since Harish-Chandra defined the Schwartz space $\mathcal{C}(G):=\mathcal{C}^{2}(G)$ and should be seen as an harvest of known results. Our approach and attainment of the full Harish-Chandra transform in Theorem $3.8$ also gives a fresh impetus to the practice of harmonic analysis in the tradition of Harish-Chandra. For example and with $rank(G)=rank(K),$ the image $\bar{\mathcal{Z}}({\mathfrak{F}}^{\epsilon})$ is known to be decomposable into a discrete part $\bar{\mathcal{Z}}_{B}({\mathfrak{F}}^{\epsilon})$ and principal part $\bar{\mathcal{Z}}_{H}({\mathfrak{F}}^{\epsilon}).$ That is, $$\bar{\mathcal{Z}}({\mathfrak{F}}^{\epsilon})
=\bar{\mathcal{Z}}_{B}({\mathfrak{F}}^{\epsilon})\times\bar{\mathcal{Z}}_{H}({\mathfrak{F}}^{\epsilon})$$ (See $[10c.],\;p.\;109,$ $[11.]$ and $[12b.],\;p.\;272$ (Theorem $53$)), where $B$ and $H$ are \textit{compact} and \textit{non-compact Cartan subgroups} of $G.$ Thus, from Corollary $3.9,$ members $\mathfrak{G}$ of the Schwartz algebra $\mathcal{C}^{p}(\widehat{G})$ are pairs $\mathfrak{G}=(\mathfrak{G}_{B},\mathfrak{G}_{H})$ given as $\mathfrak{G}_{B}=(\mathcal{H}\xi_{1})^{-1}\cdot h_{B}\cdot(\mathcal{H}\xi_{1})^{-1}$ and $\mathfrak{G}_{H}=(\mathcal{H}\xi_{1})^{-1}\cdot h_{H}\cdot(\mathcal{H}\xi_{1})^{-1}$ (with $h_{B}\in\bar{\mathcal{Z}}_{B}({\mathfrak{F}}^{\epsilon})$ and $h_{H}\in\bar{\mathcal{Z}}_{H}({\mathfrak{F}}^{\epsilon})$) both of which are linearly related as given in $[10c.],\;p.\;109$ (Definition $3$). Thence, to answer one of the questions raised in the paragraph above, the image $\mathcal{C}^{p}(\widehat{G})$ of the non-spherical Harish-Chandra Fourier transform on $G$ has a decomposition
$\mathcal{C}^{p}(\widehat{G})\cong\mathcal{C}^{p}_{B}(\widehat{G})\times\mathcal{C}^{p}_{H}(\widehat{G})$ where
$$\mathcal{C}^{p}_{B}(\widehat{G})=\{\mathfrak{G}_{B}:\mathfrak{G}_{B}=(\mathcal{H}\xi_{1})^{-1}\cdot h_{B}\cdot(\mathcal{H}\xi_{1})^{-1},\;h_{B}\in\bar{\mathcal{Z}}_{B}({\mathfrak{F}}^{\epsilon})\}$$ and $$\mathcal{C}^{p}_{H}(\widehat{G})=\{\mathfrak{G}_{H}:\mathfrak{G}_{H}=(\mathcal{H}\xi_{1})^{-1}\cdot h_{H}\cdot(\mathcal{H}\xi_{1})^{-1},\;h_{H}\in\bar{\mathcal{Z}}_{H}({\mathfrak{F}}^{\epsilon})\}.$$

It may be noted in passing that the second-half of the proof of Theorem $3.8$ reveals that the general form of the \textit{wave-packets} $\psi_{j}$ on $G,$ corresponding to any $j=(\mathcal{H}\xi_{1})^{-1}\cdot a\cdot(\mathcal{H}\xi_{1})^{-1}\in\bar{\mathcal{Z}}(K/M\times{\mathfrak{F}}^{\epsilon}\times K/M)$ (with $a\in\bar{\mathcal{Z}}({\mathfrak{F}}^{\epsilon})$), is given as
$$\psi_{j}=\mathcal{H}^{-1}((\mathcal{H}\xi_{1})^{-1})
\ast\mathcal{H}^{-1}\left((\mathcal{H}\xi_{1})\cdot j\cdot(\mathcal{H}\xi_{1})\right)\ast\mathcal{H}^{-1}((\mathcal{H}\xi_{1})^{-1}).$$ That is (by eliminating $j$),
$$\psi_{j}=\mathcal{H}^{-1}((\mathcal{H}\xi_{1})^{-1})
\ast\mathcal{H}^{-1}a\ast\mathcal{H}^{-1}((\mathcal{H}\xi_{1})^{-1});$$ from which we have earlier seen (in the proof of Theorem $3.8$) that $$\mathcal{H}(\psi_{j})=j.$$ We observe here that the \textit{general wave-packets} $\psi_{j}$ on $G$ are expressible in terms of the (normalized) \textit{spherical wave-packets} $\mathcal{H}^{-1}a$ on $G//K$ and that the \textit{general wave-packets} $\psi_{j}$ on $G$ assumes a decomposition into a convolution of three \textit{wave-packets;} namely $\mathcal{H}^{-1}((\mathcal{H}\xi_{1})^{-1}),$ $\mathcal{H}^{-1}a$ and $\mathcal{H}^{-1}((\mathcal{H}\xi_{1})^{-1}),$ with the \textit{spherical wave-packets} $\mathcal{H}^{-1}a$ (with $a\in\bar{\mathcal{Z}}({\mathfrak{F}}^{\epsilon})$) at the center. The above split of $\mathcal{C}^{p}(\widehat{G})\cong\mathcal{C}^{p}_{B}(\widehat{G})\times\mathcal{C}^{p}_{H}(\widehat{G})$ into discrete and principal parts is equally possible for the general wave-packets $\psi_{j}$ and this also answered one of the questions raised in the first paragraph after Corollary $3.9.$ Details of these shall be considered in a forthcoming paper.

The next two Corollaries show how to recover Theorems $2.2$ and $2.3$ from Theorem $3.8$ and also gives corresponding restricted form of the decomposition in Corollary $3.9$ for Theorem $2.3.$

\textbf{Corollary 3.10.} \textit{The Harish-Chandra Fourier transform of $\mathcal{C}^{p}(G/K)$ has a decomposition into a product of a function on $K/M$ with another function on $\mathfrak{F}^{\epsilon}.$ That is, $$\bar{\mathcal{Z}}(K/M\times {\mathfrak{F}}^{\epsilon})=(\mathcal{H}\xi_{1})^{-1}\cdot\bar{\mathcal{Z}}
({\mathfrak{F}}^{\epsilon}).\;\Box$$}

\textbf{Corollary 3.11}(Trombi-Varadarajan theorem). Let $0<p\leq 2,$ then \textit{$$\mathcal{H}(\mathcal{C}^{p}(G//K))\cong\bar{\mathcal{Z}}
({\mathfrak{F}}^{\epsilon}).$$}

\textbf{Proof.} We have from Lemma $3.1$ that every $f\in\mathcal{C}^{p}(G//K)$ satisfies $\xi_{1}\ast f\ast\xi_{1}=f.$ Hence, $\mathcal{H}\xi_{1}\cdot \mathcal{H}f\cdot\mathcal{H}\xi_{1}=\mathcal{H}f;$ so that $$\mathcal{H}f=(\mathcal{H}\xi_{1})^{-1}\cdot \mathcal{H}f\cdot(\mathcal{H}\xi_{1})^{-1}\in\mathcal{C}^{p}(\widehat{G}),$$ by Theorem $3.8.$ As $\mathcal{C}^{p}(G//K)$ is a proper closed subalgebra of $\mathcal{C}^{p}(G),$ the continuity of $\mathcal{H}$ ($[6.],\;p.\;340$) means that $\mathcal{H}f$ resides in a closed subalgebra of $\mathcal{C}^{p}(\widehat{G}).$ Now combining the fact that $\mathcal{H}f\in\bar{\mathcal{Z}}({\mathfrak{F}}^{\epsilon})$ ($[6.],$ Theorem $7.8.6$) with Corollary $3.7$ in the presence of Theorem $3.8$ gives the result$.\;\Box$

Analogous argument to Corollary $3.11$ may also be given to prove the main result of $[3b.].$ The next result may be seen from $[1c],\;p.\;17,$ $[12a.],\;p.\;364$ and $[3b.],\;p.\;193.$

\textbf{Lemma 3.12.} \textit{The Harish-Chandra Fourier transform of every function $f\in\mathcal{C}^{p}(G)$ may be computed (for $k_{1},k_{2}\in K,\;\lambda\in\mathfrak{F}_{I}$) as $$(\mathcal{H}f)(k_{1}M:\lambda:k_{2}M)=(\mathcal{H}\xi_{1})^{-1}\cdot\int_{AN}f(an)e^{(-\lambda+\rho)(\log a)}dadn\cdot(\mathcal{H}\xi_{1})^{-1}.$$ In particular, $\mathcal{H}f$ is independent of $k_{1},k_{2}\in K.\;\Box$}

\textbf{Corollary 3.13.} \textit{For $k_{1},k_{2}\in K,\;\lambda\in\mathfrak{F}_{I}$ we have $$(\mathcal{H}\xi_{1})(k_{1}M:\lambda:k_{2}M)=\int_{AN}e^{-2\lambda(\log a)}dadn.$$}

\textbf{Proof.} Compute using the espression $$(\mathcal{H}f)(k_{1}M:\lambda:k_{2}M)=\int_{AN}f(an)e^{(-\lambda+\rho)(\log a)}dadn.\;\Box$$

The setting for the harmonic analysis of $\mathcal{C}^{p}(G:F)$ in $[10c.]$ is that of the space $$\{(\mathcal{H}\xi_{F})^{-1}\cdot h\cdot(\mathcal{H}\xi_{F})^{-1}:\;h\in\bar{\mathcal{Z}}({\mathfrak{F}}^{\epsilon})\},$$ for any finite $F\subset\mathcal{E}(K)$ and with the restriction on the $K-$type. We however note that, since this restriction on $K-$type is a slight generalization of and reduces to the condition of $K-$biinvariance of a spherical function, it follows that Trombi's spaces $\mathcal{C}^{p}(G:F)$ and $\mathcal{C}^{p}(\widehat{G}:F)$ are not that far from Trombi-Varadarajan's spaces $\mathcal{C}^{p}(G//K)$ and $\bar{\mathcal{Z}}({\mathfrak{F}}^{\epsilon})\;([10b.],\;p.\;291),$ respectively.

Now let $\mathcal{C}^{p}(G)'$ and $\mathcal{C}^{p}(\widehat{G})'$ denote the respective topological dual spaces of $\mathcal{C}^{p}(G)$ and $\mathcal{C}^{p}(\widehat{G}),$ which are topological vector spaces in the weak topology ($[3b.],\;p.\;214$ and $[3c.]$). A distribution on $G$ will be said to be $p-$\textit{tempered} if it extends to a continuous linear functional from $\mathcal{C}^{p}(G)$ to $\C.$ A $2-$tempered distribution is simply called tempered. The precise meaning for the Harish-Chandra Fourier transform of a $p-$tempered distribution on $G$ is immediate from the following.

\textbf{Theorem 3.14.} \textit{The transpose $\mathcal{H}':\mathcal{C}^{p}(\widehat{G})'\rightarrow\mathcal{C}^{p}(G)'$ of $\mathcal{H}$ is a linear topological isomorphism$.\;\Box$}

The full \textit{invariant} harmonic analysis on $G$ is now attainable; See \textit{Coda} in $[12b.]$ and $[10d.].$ It is to be noted here that the results of $[10d.]$ depends on and is therefore restricted by those results in $[10c.],$ where the Harish-Chandra Fourier transform $f\mapsto\mathcal{H}f$ was considered only for $f\in\mathcal{C}^{p}(G:F)$ ($\subset\mathcal{C}^{p}(G)).$ We shall however show here that in the presence of the above Theorem $3.8$ (which is valid for all $f\in\mathcal{C}^{p}(G)$) a full invariant harmonic analysis on $G$ may now be developed as follows, with proof essentially as in $[10d.].$

To this end, let $\theta_{\pi}$ denote the global characters of a quasi-simple admissible representation $\pi$ of $G$ and, for $f\in\mathcal{C}^{p}(G),$ write $\widehat{f}$ which is defined on $\widehat{G}$ as $$\widehat{f}(\pi)=\theta_{\pi}(f)=\int_{G'}f(x^{-1})\theta_{\pi}(x)dx$$ ($G'$ being the regular set in $G$ on which $\theta_{\pi}$ is well-known to be \textit{analytic}) and is termed the invariant Harish-Chandra Fourier transform of $f.$ Since the global character $\theta_{\pi}$ is the distribution on $G$ given as $\theta_{\pi}(f)=tr(\int_{G}f(x)\pi(x)dx)$ we then have that $$\widehat{f}(\pi)=tr(\int_{G'}f(x)\pi(x)dx)=tr(\mathcal{H}f),\;f\in\mathcal{C}^{p}(G).$$ It has been shown by Trombi $[8'c.,d.]$ that the split into discrete and principal parts is to be expected even at the level of $\mathcal{C}^{p}(\widehat{G})$ (See also Theorems $60$ and $65$ in $[12b.]$). We therefore proceed as follows.

We denote the principal and discrete series of representations of $G$ by $\pi_{\sigma,\lambda}$ (with $\sigma\in\widehat{M},$ $\lambda\in int(\mathfrak{F}^{\epsilon})$) and $\pi_{\omega}$ (with $\omega\in\widehat{G}_{d}$), respectively. The following results are well-known.

\textbf{Lemma 3.15.} ($[10d.]$) \textit{Let $0<p\leq 2$ and $f\in\mathcal{C}^{p}(G).$ Then\\
$(i.)$ $\widehat{f}(\sigma:\cdot)$ is an entire function on $\mathfrak{F}_{c}$ of exponential type;\\
$(ii.)$ $\widehat{f}(s\sigma:s\lambda)=\widehat{f}(\sigma:\lambda)$ for all $s\in\mathfrak{w},\;(\sigma,\lambda)\in\widehat{M}\times\mathfrak{F}_{c};$ \\
$(iii.)$ $\widehat{f}(\sigma:\lambda)=0,$ if $\sigma\notin\widehat{M}$ and
$\widehat{f}(\omega)=0,$ if $\omega\notin\widehat{G}_{d}.$}

\textbf{Proof.} It is sufficient to prove these results for $f\in C^{\infty}_{c}(G).$ See $[10d.]$ and $(3.3.9)$ of $[6.].\;\Box$

The properties included in Lemma $3.15$ suggest a candidate for the image of the invariant Harish-Chandra Fourier transform. Set $$C^{p}(G)=\{\theta_{\sigma,\lambda}:\sigma\in\widehat{M},\lambda\in int(\mathfrak{F}^{\epsilon}_{c})\}\bigcup\{\theta_{\omega}:\omega\in\widehat{G}_{d}\}\bigcup\mathfrak{B}_{p}$$ where $\mathfrak{B}_{p}$ is as in $[10d.]$ (See also $[12b.],\;p.\;285$).

\textbf{Definition 3.16} ($[10d.]$) \textit{Let $\mathcal{C}^{p}(C(G))_{o}$ denote the linear space of all complex-valued functions on $C^{p}(G)$ such that for $L\in\mathcal{C}^{p}(C(G))_{o}$ (and by denoting $L(\theta_{\sigma,\lambda})$ as $L(\sigma,\lambda)$ and $L(\theta_{\omega})$ as $L(\omega)$);\\
$(i.)$ each $L(\sigma:\cdot)$ is holomorphic on $int(\mathfrak{F}^{\epsilon}_{c})$;\\
$(ii.)$ $L(s\sigma:s\lambda)=L(\sigma:\lambda)$ for all $s\in\mathfrak{w},\;(\sigma,\lambda)\in\widehat{M}\times int(\mathfrak{F}^{\epsilon}_{c});$ \\
$(iii.)$ $L(\sigma:\lambda)=0,$ if $\sigma\notin\widehat{M}$ and
$L(\omega)=0,$ if $\omega\notin\widehat{G}_{d};$\\
$(iv.)$ $\sup_{\widehat{M}\times int(\mathfrak{F}^{\epsilon}_{c})}
(1+\mid\lambda \mid)^{\alpha}\mid L(\sigma:\lambda;u)\mid=:\nu^{p}_{u,\alpha}(L)<\infty,$ for all $\alpha\in\R,$ $u\in S(\mathfrak{F}_{c}).$}

Now let $\mathcal{C}^{p}(C(G))$ denote the subspace of functions $L\in\mathcal{C}^{p}(C(G))_{o}$ such that $$L(\sigma;t\zeta;\partial^{k}(u))=\sum_{\theta\in\mathfrak{B}_{p}}c_{p}(\theta_{\sigma,t\zeta,k}:\theta)L(\theta),$$ where
$\sigma\in\widehat{M},\;t\in\mathfrak{w},\;\zeta\in V_{p},\;0\leq k\leq\theta_{t}(\zeta)-1,\;u\in S(\mathfrak{F}_{c}).$ The reader is referred to the remarks on $p.\;285$ of $[12b.]$ following the definition of $\widehat{f}$ for the motivation of and necessity for the set $\mathfrak{B}_{p}$ and the requirement on $L$ given above in the case of $SL(2,\R)$. We then give $\mathcal{C}^{p}(C(G))$ the topology generated by the family of seminorms $\mu^{p}_{u,\alpha}$  given as $$\mu^{p}_{u,\alpha}(L)=\nu^{p}_{u,\alpha}(L)+(\sum_{\omega\in \widehat{G}_{d}}\mid L(\omega)\mid^{2})^{\frac{1}{2}},$$ $u\in S(\mathfrak{F}_{c}),\;\alpha\in\R.$

\textbf{Lemma 3.17.} \textit{The family $\mu^{p}_{u,\alpha}$ of seminorms convert $\mathcal{C}^{p}(C(G))$ into a Schwartz space$.\;\Box$}

\textbf{Theorem 3.18.} \textit{The invariant Harish-Chandra Fourier transform $f\mapsto\widehat{f}$ is a linear topological algebra isomorphism $\mathcal{C}^{p}(G)\rightarrow\mathcal{C}^{p}(C(G)).$}

\textbf{Proof.} The map $tr$ is clearly continuous, using an argument analogous to Proposition $1$ of $[10d.],\;p.\;235,$ while the surjectivity argument for $f\mapsto\widehat{f}$ in Theorem $1$ of $[10d.],\;p.\;235$ also holds for all $\mathcal{C}^{p}(G),$ now that we already have the \textit{full} isomorphism $\mathcal{H}:\mathcal{C}^{p}(G)\cong\mathcal{C}^{p}(\widehat{G})$ in Theorem $3.8.\;\Box$

We now have the linear toplogical algebra isomorphisms $$\mathcal{C}^{p}(G)\cong\mathcal{C}^{p}(\widehat{G})\cong\mathcal{C}^{p}(C(G))\;\mbox{given by}\;f\mapsto\mathcal{H}f\mapsto\widehat{f}.$$

\textbf{Proposition 3.19.} \textit{Members of $\mathcal{C}^{p}(C(G))$ are all of the form
$$L((\sigma;\lambda)\oplus\omega)=tr\left[((\mathcal{H}\xi_{1})^{-1})^{2}\cdot\int_{AN}f(an)e^{(\lambda+\rho)(\log a)}dadn\right],$$ for any $f\in\mathcal{C}^{p}(G).$ More precisely, we have that $L(\sigma,\lambda)=tr\left[((\mathcal{H}\xi_{1})^{-1})^{2}\cdot h_{H}\right],$ with $h_{H}\in\bar{\mathcal{Z}}_{H}({\mathfrak{F}}^{\epsilon})$ and $L(\omega)=tr\left[((\mathcal{H}\xi_{1})^{-1})^{2}\cdot h_{B}\right],$ with $h_{B}\in\bar{\mathcal{Z}}_{B}({\mathfrak{F}}^{\epsilon}).$}

\textbf{Proof.} We combine the remarks following Corollary $3.9$ with Lemma $3.12.\;\Box$

We shall refer to a map $m$ of functions on $G$ as being \textit{invariant} whenever $m(f^{y})=m(f),$ where $f^{y}(x)=f(y^{-1}xy),\;x,y\in G.$ The following result is the reason for the term \textit{invariant harmonic analysis.}

\textbf{Theorem 3.20.} \textit{The invariant Harish-Chandra Fourier transform $f\mapsto\widehat{f}$ is invariant in the above sense.}

\textbf{Proof.} We set $y_{1}=y$ and $y_{2}=y^{-1}$ in the formula on $p.\;(1.6)$ of $[1c.],$ to get
$$(\mathcal{H}f^{y})(\lambda)=\pi_{\sigma,\lambda}(y)\cdot\mathcal{H}f\cdot\pi_{\sigma,\lambda}(y^{-1}),$$ for all $(\sigma,\lambda)\in\widehat{M}\times int(\mathfrak{F}^{\epsilon}),\;f\in\mathcal{C}^{p}(G).$ Hence, we have $\widehat{f^{y}}=tr(\mathcal{H}f^{y})
=tr(\pi_{\sigma,\lambda}(y)\cdot\mathcal{H}f\cdot\pi_{\sigma,\lambda}(y^{-1}))
=tr(\pi_{\sigma,\lambda}(y)\cdot\pi_{\sigma,\lambda}(y^{-1})\cdot\mathcal{H}f)=tr(\mathcal{H}f)=\widehat{f}.\;\Box$

Denote the kernel of the map $tr$ in $\mathcal{C}^{p}(G)$ by $K^{p}(G).$ It is known that $K^{2}(SL(2,\R))$ is the closure of the span of the commutators in $\mathcal{C}^{2}(SL(2,\R));$ $[12b.],$ Theorem $65,p.\;289.$ Then $\mathcal{C}^{p}(C(G))/K^{p}(G)$ is a commutative Fr$\acute{e}$chet algebra with operation induced from the convolution on $\mathcal{C}^{p}(C(G)).$ The duo of Theorem $3.18$ and and the Fundamental Theorem of homomorphisms imply that $$\mathcal{C}^{p}(G)/K^{p}(G)\cong \mathcal{C}^{p}(C(G)).$$ See also Theorem $66$ of $[12b.].$ The following is a generalization of \textit{Rao-Varadarajan theorem} on $\mathcal{C}^{2}(SL(2,\R))$ to all of $\mathcal{C}^{2}(G)$.

\textbf{Theorem 3.21.} \textit{$f\in K^{p}(G)$ if and only if $\theta(f)=0,$ for every $\theta\in C^{p}(G).$ A $p-$tempered distribution on $G$ is invariant if and only if it vanishes on $K^{p}(G).$}

\textbf{Proof.} Our argument is as in $[12b.],p.\;289.\;\Box$

A further realization of members of $\mathcal{C}^{p}(C(G))$ than in Proposition $3.19$ and an explicit computation of the members of $K^{p}(G)$ are desirable. A way forward is in the computation of $(\mathcal{H}\xi_{1})^{-1}.$ Let us now consider the zero-Schwartz spaces, $\mathcal{C}^{0}(G)$ and $\mathcal{C}^{0}(\widehat{G}).$

We define $$\mathcal{C}^{0}(G)=\bigcap_{0<p\leq2}\mathcal{C}^{p}(G)$$ from which it follows that $$C^{\infty}_{c}(G)\subset\mathcal{C}^{0}(G)
\subset\cdots\subset\mathcal{C}^{p_{1}}(G)\subseteq\mathcal{C}^{p_{2}}(G)\subset\cdots$$ with $0<p_{1}\leq p_{2}\leq2.$ We topologize $\mathcal{C}^{0}(G)$ with the \textit{projective limit topology} for this intersection. Each of these subspaces is a Fr$\acute{e}$chet algebra under convolution (See $[2.]$ and the references contained in section $19$), $C^{\infty}_{c}$ is dense in $\mathcal{C}^{0}(G)$ with continuous inclusion. This shows that (both of) the (invariant and non-invariant) Harish-Chandra Fourier transform $\mathcal{H}$ on $\mathcal{C}^{p}(G)$ may be restricted to $\mathcal{C}^{0}(G).$ We also define $\mathcal{C}^{0}(\widehat{G})$ in a completely analogous manner as done for $\mathcal{C}^{0}(G).$ Hence, $\mathcal{H}$ of Theorem $3.8$ restricts to $\mathcal{C}^{0}(G)$ and we have $\mathcal{H}:\mathcal{C}^{0}(G)\rightarrow\mathcal{C}^{0}(\widehat{G}).$

The spaces $\mathcal{C}^{0}(G)$ and $\mathcal{C}^{0}(\widehat{G})$ may also be topologized by means of other seminorms, instead of the projective limit topologies for their corresponding intersections ($[2.],\;p.\;99$ and $p.\;102$). We however have the following result in any of the said equivalent topologies.

\textbf{Theorem 3.22.} \textit{The Harish-Chandra Fourier transform $\mathcal{H}$ sets up a linear topological algebra isomorphism $\mathcal{H}:\mathcal{C}^{0}(G)\rightarrow\mathcal{C}^{0}(\widehat{G}).$}

\textbf{Proof.} We simply take $\mathcal{H}:\mathcal{C}^{0}(G)\rightarrow\mathcal{C}^{0}(\widehat{G})$ as the restriction of the linear topological algebra isomorphism in Theorem $3.8.\;\Box$

\ \\
{\bf \S4. The example of spherical convolutions; $g_{\lambda,A}:=\widetilde{f}\ast\varphi_{\lambda}.$}

Let $f\in\mathcal{C}^{p}(G)$ and define $\widetilde{f}$ by $\widetilde{f}=f_{\mid_{A^{+}}},$ then $\widetilde{f}\in\mathcal{C}^{p}(A^{+}).$

\textbf{Lemma 4.1.} \textit{Set $\mathcal{A}(U:\chi)=\{f\in C^{\infty}(U):\delta'(q)f=\chi(q)f,\;\forall\;q\in\mathfrak{D}\}$ for any open set $U$ in $A^{+}$ and homomorphism $\chi:\mathfrak{D}\rightarrow\C.$ Then we have $g_{\lambda,A}\in \mathcal{A}(U:\gamma(\cdot)(\lambda)).$}

\textbf{Proof.} We know that if $f,g\in C^{\infty}(U)$ and $a,b\in\mathfrak{U}(\mathfrak{g}_{\C})$ then $$a(f\ast g)b=fb\ast ag,$$ $[6.],\;p.\;255.$ Hence, $$\delta'(q)g_{\lambda,A}=\widetilde{f}\ast (\delta'(q)\varphi_{\lambda})
=\gamma(q)(\lambda)(\widetilde{f}\ast\varphi_{\lambda})=\gamma(q)(\lambda)(g_{\lambda,A}).\;\Box$$

This Lemma shows that $g_{\lambda,A}$ is an eigenfunction on $A^{+}.$ Hence, $g_{\lambda,A}$ has a Harish-Chandra series expansion on $A^{+}$ which is $$g_{\lambda,A}(h)=\sum_{s\in\mathfrak{w}}c(s\lambda)\left(e^{(s\lambda-\rho)(\log h)}+\sum_{\mu\in L^{+}}a_{\mu}(s\lambda)e^{(s\lambda-\rho-\mu)(\log h)}\right),$$ valid for all $h\in A^{+},$ and some $\lambda,$ $L^{+}$ in $[6.],$ with the $c-$function now given as $c(s\lambda)=\sum_{1\leq i\leq w}\gamma^{si}(\lambda)g_{\lambda,A}(h_{o};u'_{i})$ where $w=\mid \mathfrak{w}\mid,s\in\mathfrak{w}$ in which $\gamma^{si}(\lambda)$ are the entries of the inverse matrix of the invertible $w\times w$ matrix $(\gamma_{si}(\lambda))_{1\leq i\leq w}$ (with $\gamma_{si}(\lambda)=\Phi(s\lambda:h_{o};u'_{i}),$ a basis for $\mathcal{A}(A^{+}:\chi_{\lambda})$), $h_{o}\in A^{+}$ and each $u'_{i}$ is the radial component of each $u_{i}\in S(\mathfrak{F}).$

In order to then compute $\mathcal{H}(\mathcal{C}^{p}(A^{+}))$ we shall employ the methods of \S $3.$ to prove the following.

\textbf{Lemma 4.2.} \textit{$\mathcal{H}\widetilde{f}=(\mathcal{H}g_{\lambda,A})\cdot(\mathcal{H}\varphi_{\lambda})^{-1}$ for $\lambda\in\mathfrak{F}_{I}.$}

\textbf{Proof.} $\mathcal{H}g_{\lambda,A}=\mathcal{H}(\widetilde{f}\ast\varphi_{\lambda})
=(\mathcal{H}\widetilde{f})\cdot(\mathcal{H}\varphi_{\lambda}).\;\Box$

This reveals that the non-spherical Harish-Chandra Fourier of $\widetilde{f}$ is given in terms of $\mathcal{H}g_{\lambda,A}$ and $(\mathcal{H}\varphi_{\lambda})^{-1}$ in which both $g_{\lambda,A}$ and $\varphi_{\lambda}$ are elementary spherical functions.

\textbf{Proposition 4.3.} \textit{$\mathcal{H}(\mathcal{C}^{p}(A^{+}))
=\{(\mathcal{H}g_{\lambda,A})\cdot(\mathcal{H}\varphi_{\lambda})^{-1}:\lambda\in\mathfrak{F}_{I}\}.\;\Box$}

Explicit computations for both $\mathcal{H}g_{\lambda,A}$ and $(\mathcal{H}\varphi_{\lambda})^{-1}$ via the consideration of their Harish-Chandra series expansions give concrete members of $\mathcal{H}(\mathcal{C}^{p}(A^{+})).$ Here $\mathcal{H}$ is the spherical Harish-Chandra Fourier transform given as $$(\mathcal{H}f)(\nu)=\int_{G}f(x)\varphi_{-\nu}(x)dx=\int_{K\cdot cl(A^{+})\cdot K}f(ka_{t}k)\varphi_{-\nu}(a_{t})J(a_{t})dkdtdk,$$ with $J(a_{t})=\prod_{\lambda\in\Delta^{+}}\left[\sinh \lambda(H)\right]^{\dim \mathfrak{g}_{\lambda}}$ ($[6.],\;p.\;73$). This is however straightforward.

We believe that the explicit computation of the function $(\mathcal{H}\xi_{1})^{-1}(\lambda),$ with $\lambda\in\mathfrak{F}^{\epsilon},$ is necessary in order to pave way for further research along our perspective. This will however be taken up in another paper. Presently, we have that $$(\mathcal{H}\xi_{1})(\lambda)=\int_{AN}e^{(-\lambda+\rho)(\log a)}\varphi_{\lambda}(an)dadn=$$ $$\int_{AN}e^{(-\lambda+\rho)(\log a)}\left[\sum_{s\in\mathfrak{w}}c(s\lambda)\left(e^{(s\lambda-\rho)(\log (an))}+\sum_{\mu\in L^{+}}a_{\mu}(s\lambda)e^{(s\lambda-\rho-\mu)(\log (an))}\right)\right]dadn$$
$$=\sum_{s\in\mathfrak{w}}c(s\lambda)\left(\int_{AN}e^{(s\lambda-\lambda)(\log (an))}dadn+\sum_{\mu\in L^{+}}a_{\mu}(s\lambda)\int_{AN}e^{(s\lambda-\lambda-\mu)(\log (an))}dadn\right)$$ (since $\log (an)=\log (a)$) whose inverse is required in Theorem $3.8.$

\ \\
\ \\
\ \\
\ \\
{\bf   References.}
\begin{description}
\item [{[1.]}] Arthur, J. G., $(a.)$ \textit{Harmonic analysis of tempered distributions on semi-simple Lie groups of real rank one}, Ph.D. Dissertation, Yale University, $1970;$ $(b.)$ \textit{Harmonic analysis of the Schwartz space of a reductive Lie group I,} mimeographed note, Yale University, Mathematics Department, New Haven, Conn; $(c.)$ \textit{Harmonic analysis of the Schwartz space of a reductive Lie group II,} mimeographed note, Yale University, Mathematics Department, New Haven, Conn.
        \item [{[2.]}] Barker, W. H.,  $L^p$ harmonic analysis on $SL(2, \R),$ \textit{Memoirs of American Mathematical Society,} $vol.$ {\bf 76} , no.: {\bf 393}. $1988$
            \item [{[3.]}] Eguchi, M., $(a.)$ The Fourier Transform of the Schwartz space on a semisimple Lie group, \textit{Hiroshima Math. J.}, $vol.$ \textbf{4}, ($1974$), pp. $133-209.$ $(b.)$ Asymptotic expansions of Eisenstein integrals and Fourier transforms on symmetric spaces, \textit{J. Funct. Anal.} \textbf{34}, ($1979$), pp. $167 - 216.$ $(c.)$ Some properties of Fourier transform on Riemannian symmetric spaces, \textit{Lecture on Harmonic Analysis on Lie Groups and related Topics,} (T. Hirai and G. Schiffmann (eds.)) Lectures in Mathematics, Kyoto University, No. \textbf{4}) pp. $9 - 43.$
                \item [{[4.]}] Eguchi, M. and Kowata, A., On the Fourier transform of rapidly decreasing function of $L^{p}$ type on a symmetric space, \textit{Hiroshima Math. J.} $vol.$ \textbf{6}, ($1976$), pp. $143 - 158.$
                \item [{[5.]}] Ehrenpreis, L. and Mautner, F. I., Some properties of the Fourier transform on semisimple Lie groups, I,
                \textit{Ann. Math.}, $vol.$ $\textbf{61}$ ($1955$), pp. 406-439; II, \textit{Trans. Amer. Math. Soc.}, $vol.$ $\textbf{84}$ ($1957$), pp. $1-55;$ III, \textit{Trans. Amer. Math. Soc.}, $vol.$ $\textbf{90}$ ($1959$), pp. $431-484.$
                    \item [{[6.]}] Gangolli, R. and Varadarajan, V. S., \textit{Harmonic analysis of spherical functions on real reductive groups,} Ergebnisse der Mathematik und iher Genzgebiete, $vol.$ {\bf 101}, Springer-Verlag, Berlin-Heidelberg. $1988.$
                    \item [{[7.]}] Helgason, S., $(a)$  A duality for symmetric spaces with applications to group representations, \textit{Advances in Mathematics,} $vol.$ $\textbf{5}$ ($1970$), pp. $1-154.$ $(b)$ \textit{Differential geometry and symmetric spaces,} Academic Press, New York, $1962.$
                        \item [{[8.]}] Knapp, A.W., \textit{Representation theory of semisimple groups; An overview based on examples,} Princeton University Press, Princeton, New Jersey. $1986.$
                \item [{[9.]}] Oyadare, O. O., $(a.)$ On harmonic analysis of spherical convolutions on semisimple Lie groups, \textit{Theoretical Mathematics and Applications}, $vol.$ $\textbf{5},$ no.: {\bf 3}. ($2015$), pp. $19-36.$ $(b.)$ Series analysis and Schwartz algebras of spherical convolutions on semisimple Lie groups (\textit{under review}), arXiv.$1706.09045$ [math.RT]
                    \item [{[10.]}] Trombi, P. C., $(a.)$ Spherical transforms on symmetric spaces of rank one (or Fourier analysis on semisimple Lie groups of split rank one), \textit{Thesis, University of Illinios} ($1970$). $(b.)$ On Harish-Chandra's theory of the Eisenstein integral for real semisimple Lie groups. \textit{University of Chicago Lecture Notes in Representation Theory,} ($1978$), pp. $287$-$350.$ $(c.)$ Harmonic analysis of $\mathcal{C}^{p}(G:F)\;(1\leq p<2)$ \textit{J. Funct. Anal.,} $vol.$ {\bf 40}. ($1981$), pp. $84$-$125.$ $(d.)$ Invariant harmonic analysis on split rank one groups with applications. \textit{Pacific J. Math.} $vol.$ \textbf{101}. no.: \textbf{1}.($1982$), pp. $223 - 245.$
                \item [{[11.]}] Trombi, P. C. and Varadarajan, V. S., Spherical transforms on semisimple Lie groups, \textit{Ann. Math.,} $vol.$ {\bf 94}. ($1971$), pp. $246$-$303.$
                        \item [{[12.]}] Varadarajan, V. S., $(a.)$ Eigenfunction expansions on semisimple Lie groups, in \textit{Harmonic Analysis and Group Representation}, (A. Fig$\grave{a}$  Talamanca (ed.)) (Lectures given at the $1980$ Summer School of the \textit{Centro Internazionale Matematico Estivo (CIME)} Cortona (Arezzo), Italy, June $24$ - July $9.$ vol. \textbf{82}) Springer-Verlag, Berlin-Heidelberg. $2010,$ pp. $351-422.$ $(b.)$ \textit{An introduction to harmonic analysis on semisimple Lie groups,} Cambridge Studies in Advanced Mathematics, \textbf{161},  Cambridge University Press, $1989.$ $(c.)$ Harmonic analysis on real reductive reductive groups, \textit{Lecture Notes in Mathematics,} \textbf{576}, Springer Verlag, $1977.$
                        \end{description}
\end{document}